\input amstex\documentstyle{amsppt}  
\pagewidth{12.5cm}\pageheight{19cm}\magnification\magstep1
\topmatter
\title On the character of certain tilting modules\endtitle
\author G. Lusztig and G. Williamson\endauthor
\address{Department of Mathematics, M.I.T., Cambridge, MA and Max Planck Institut f\"ur Mathematik, Bonn}
\endaddress
\thanks{G.L. supported in part by National Science Foundation grant DMS-1303060 and by a Simons Fellowship.}
\endthanks
\endtopmatter   
\document
\define\Rep{\text{\rm Rep}}

\define\frl{\forall}

\define\qua{\quad}

\define\part{\partial}

\define\m{\mapsto}
\define\do{\dots}

\define\lra{\leftrightarrow}

\define\sub{\subset}

\define\Hom{\text{\rm Hom}}

\define\Aut{\text{\rm Aut}}

\define\a{\alpha}

\define\e{\epsilon}

\define\r{\rho}

\redefine\l{\lambda}
\define\z{\zeta}
\define\x{\xi}

\define\kk{\bold k}

\define\NN{\bold N}

\define\ZZ{\bold Z}

\define\ct{\Cal T}

\define\che{\check}

\define\AN{A1}
\define\ANI{A2}
\define\DO{Do}
\define\JE{Je}
\define\KL{KL}
\define\LU{Lu}
\define\SOI{S1}
\define\SOII{S2}
\define\SOIII{S3}

\subhead 1 \endsubhead
Let $G$ be an almost simple, simply connected algebraic group over $\kk$, an algebraically closed field of 
characteristic $p>1$. We fix a maximal torus $T$ of $G$; let $X=\Hom(T,\kk^*)$ (with group operation 
written as $+$). Let $\Rep G$ be the category of finite dimensional $\kk$-vector spaces with a given 
rational linear action of $G$ and let $\ct G$ be the full subcategory of $\Rep G$ consisting of tilting
modules (see \cite{\DO}). Let $X^+$ be the set of dominant elements in $X$, defined in terms of a fixed
choice of a set $R^+$ of positive roots in the set of roots of $G$ with respect to $T$. For $V\in\Rep G$ and
$\mu\in X$ let $n_\mu(V)$ be the dimension of the $\mu$-weight space of $V$. For $V\in\Rep G$ we set
$[V]=\sum_{\mu\in X}n_\mu(V)e^\mu\in\ZZ[X]$, where $\ZZ[X]$ is the group ring of $X$ with standard basis
$\{e^\mu;\mu\in X\}$. (Actually, we have $[V]\in\ZZ[X]^W$ where $W$ is the Weyl group of $G$ with respect to
$T$, viewed as a subgroup of $\Aut(X)$ and $\ZZ[X]^W$ denotes the subring of $W$-invariants.) According to 
\cite{\DO}, there is a well defined bijection $\l\lra T_\l$ between $X^+$
and a set of representatives for the isomorphism classes of indecomposable objects of $\ct G$ such that the 
following holds: $n_\l(T_\l)=1$; moreover, if $\mu\in X$, $n_\mu(T_\l)\ne0$ then $\mu\le\l$ ($\le$ as in
no.2). 
We set $X^{++}=X^++\r$ where $\r\in X^+$ is defined by 
$2\r=\sum_{\a\in R^+}\a$. For $\z\in X^{++}$ we set $S_\z=[T_{\z-\r}]$,
$$S_\z^0=\sum_{w\in W}\e_we^{w(\z)}/\sum_{w\in W}\e_we^{w(\r)}\in\ZZ[X]^W.$$
(Here $w\m\e_w$ is the sign character of $W$.) By general principles, for $\z\in X^{++}$ we have
$$S_\z=\sum_{\mu\in X^{++}}y_{\mu,\z}S_\mu^0$$
where $y_{\mu,\z}\in\NN$ is zero unless $\mu\le\z$. 
One of the mysteries of the subject is that there is no known explicit formula (even conjecturally and even 
for $p\gg0$) for the coefficients $y_{\mu,\z}$ above. This in contrast with the situation for the quantum 
group at ${}^p\sqrt{1}$ when, for $\l\in X^+$, the analogue $T'_\l$ of $T_\l$ and the corresponding element 
$[T'_\l]\in\ZZ[X]^W$ which records weight multiplicities in $T'_\l$ are defined (see \cite{\AN}) and, 
setting $S^1_\z=[T'_{\z-\r}]$ for $\z\in X^{++}$, we have
$$S^1_\z=\sum_{\mu\in X^{++}}x_{\mu,\z}S_\mu^0\tag a$$
where $x_{\mu,\z}\in\NN$ can be explicitly computed in terms of the 
polynomials $Q_{y,w}$ of \cite{\KL}, see Soergel \cite{\SOI,\SOII,\SOIII} (it is zero unless $\mu\le\z$).

Now \cite{\LU} gives a method to express directly the characters of simple objects of $\Rep G$ (for $p\gg0$)
in terms of the characters of the irreducible representation of the quantum group, without the use of the 
Steinberg tensor product theorem. In this paper we show that a variation of the method of \cite{\LU}, 
combined with results of Donkin and Andersen, gives (at least conjecturally and for $p\gg0$), a simple 
closed formula for the coefficients $y_{\mu,\z}$ (for a large set of $\z\in X^{++}$) in terms of the 
coefficients $x_{\mu,\z}$. 

We thank H. H. Andersen for help with references.

\subhead 2. Notation \endsubhead
For $\l,\mu\in X$ we write $\mu\le\l$ if $\l-\mu\in\sum_{\a\in R^+}\NN\a$.

Let $\{\che\a_i;i\in I\}$ be the basis of $\Hom(X,\ZZ)$ consisting of simple coroots. We have 
$X^+=\{\l\in X;\che\a_i(\l)\in\NN\qua\frl i\in I\}$.
For any $\x=\sum_{\l\in X}c_\l e^\l\in\ZZ[X]$ (with $c_\l\in\ZZ$) and any $h\ge0$ we set
$\x^{(h)}=\sum_{\l\in X}c_\l e^{p^h\l}\in\ZZ[X]$. Let 
$$X^+_{red}=\{\l\in X;\che\a_i(\l)\in[0,p-1]\qua\frl i\in I\}.$$
For any $\l\in X^+$ we define $\l^0,\l^1,\l^2,\do$ in $X^+_{red}$ by $\l=\sum_{k\ge0}p^k\l^k$; note that
$\l^k=0$ for large $k$. 

\subhead 3 \endsubhead
For $\l\in X^+,\nu\in X^+_{red},\mu\in X^+$ such that $\l+\r=\nu+p(\mu+\r)$ we have the Donkin tensor
product formula, see \cite{\DO}: 
$$[T_\l]=[T_{\nu+(p-1)\r}][T_\mu]^{(1)},\tag a$$
and its quantum analogue, see \cite{\AN}:
$$[T'_\l]=[T'_{\nu+(p-1)\r}](S^0_{\mu+\r})^{(1)}.\tag b$$
For $k\ge0$ let $X^+_k=\{\z\in X^+;\z^k\in X^{++},\z^{k+1}=\z^{k+2}=\do=0\}$. Note that $X^+_k\sub X^{++}$.
For any $k\ge0$ and any $\z\in X^+_k$ we show:
$$S_\z=S_{\z^0+p\r}S_{\z^1+p\r}^{(1)}\do S_{\z^{k-1}+p\r}^{(k-1)}S_{\z^k}^{(k)}.\tag c$$
We argue by induction on $k$. For $k=0$ the result is obvious: we have $z=z^0\in X^{++}$ and 
$S_\z=S_{\z^0}$. Assume now that $k\ge1$. We have $\z=\z^0+p(\z^1+p\z^2+\do)$ hence, by (a):
$$S_\z=S_{\z^0+p\r}S_{\z^1+p\z^2+\do}^{(1)}.$$
We have $\z^1+p\z^2+\do\in X^+_{k-1}$ hence by the induction hypothesis we have
$$S_{\z^1+p\z^2+\do}=S_{\z^1+p\r}\do S_{\z^{k-1}+p\r}^{(k-2)}S_{\z^k}^{(k-1)}.$$
The result follows.

Now, using (b), we see that, if $\z\in X^+$ satisfies $\z^1+p\z^2+\do\in X^{++}$ then
$$S^1_\z=S^1_{\z^0+p\r}(S^0_{\z^1+p\z^2+\do})^{(1)}.\tag d$$

\subhead 4 \endsubhead
For any $h\ge1$ and any $\z\in X^+$ such that $\z^{h-1}+p\z^h+p^2\z^{h+1}+\do\in X^{++}$ we define $S^h_\z$ 
by the inductive formula
$$S^h_\z=\sum_{\mu\in X^{++}}x_{\mu,\z^{h-1}+p\z^h+p^2\z^{h+1}+\do}
S^{h-1}_{\z^0+p\z^1+\do+p^{h-2}\z^{h-2}+p^{h-1}\mu}.\tag a$$
If $h=1$ this agrees with 0(a) (the condition on $\z$ becomes $\z\in X^{++}$).
If $h\ge2$, then $S^{h-1}_{\z^0+p\z^1+\do+p^{h-2}\z^{h-2}+p^{h-1}\mu}$ is defined by induction since
$$\z':=\z^0+p\z^1+\do+p^{h-2}\z^{h-2}+p^{h-1}\mu$$ 
satisfies 
$$\z'{}^{h-2}+p\z'{}^{h-1}+p^2\z'{}^h+\do=\z^{h-2}+p\mu\in X^{++}.$$
For any $h\ge1$ and any $\z\in X^+$ such that $\z^h+p\z^{h+1}+p^2\z^{h+2}+\do\in X^{++}$ we show:
$$S^h_\z=S^1_{\z^0+p\r}(S^1_{\z^1+p\r})^{(1)}\do(S^1_{\z^{h-1}+p\r})^{(h-1)}
(S^0_{\z^h+p\z^{h+1}+\do})^{(h)}.\tag b$$
Note that our assumption implies $\z^{h-1}+p\z^h+p^2\z^{h+1}+\do\in X^{++}$ hence $S^h_\z$ is defined. We 
argue by induction on $h$. Assume first that $h=1$ and $\z^1+p\z^2+p^2\z^3+\do\in X^{++}$. Then (b) reduces 
to
$S^1_\z=S^1_{\z^0+p\r}(S^0_{\z^1+p\z^2+\do})^{(1)}$ which is known from 3(d). Now assume that $h\ge2$. In the
right hand side of (a) we replace (using the induction hypothesis)
$S^{h-1}_{\z^0+p\z^1+\do+p^{h-2}\z^{h-2}+p^{h-1}\mu}$ by
$$S^1_{\z^0+p\r}(S^1_{\z^1+p\r})^{(1)}\do(S^1_{\z^{h-2}+p\r})^{(h-2)}(S^0_\mu)^{(h-1)}.$$
(Note that, if $\z'$ is as above, then $\z'{}^{h-1}+p\z'{}^h+\do=\mu\in X^{++}$, hence the induction 
hypothesis is applicable.) Thus from (a) we obtain
$$\align&S^h_\z\\&=\sum_{\mu\in X^{++}}x_{\mu,\z^{h-1}+p\z^h+p^2\z^{h+1}+\do}
S^1_{\z^0+p\r}(S^1_{\z^1+p\r})^{(1)}\do(S^1_{\z^{h-2}+p\r})^{(h-2)}(S^0_\mu)^{(h-1)}.\endalign$$
It remains to show
$$\sum_{\mu\in X^{++}}x_{\mu,\z^{h-1}+p\z^h+p^2\z^{h+1}+\do}S^0_\mu
=S^1_{\z^{h-1}+p\r}(S^0_{\z^h+p\z^{h+1}+\do})^{(1)}$$
that is,
$$S^1_{\z^{h-1}+p\z^h+p^2\z^{h+1}+\do}=S^1_{\z^{h-1}+p\r}(S^0_{\z^h+p\z^{h+1}+\do})^{(1)}.$$
This is known from 3(d). Thus, (a) is proved.

\subhead 5 \endsubhead
Now assume that $k\ge1$ and $\z\in X^+_k$ that is, $\z^k\in X^{++}$, $\z^{k+1}=\z^{k+2}=\do=0$. Then 
$S^{k+1}_\z$ is defined. We have
$$S^{k+1}_\z=\sum_{\mu\in X^{++}}x_{\mu,\z^k}S^k_{\z^0+p\z^1+\do+p^{k-1}\z^{k-1}+p^k\mu}.$$
By 4(b), this becomes
$$S^{k+1}_\z=\sum_{\mu\in X^{++}}x_{\mu,\z^k}S^1_{\z^0+p\r}(S^1_{\z^1+p\r})^{(1)}\do
(S^1_{\z^{k-1}+p\r})^{(k-1)}(S^0_\mu)^{(k)}.$$
Hence, using 3(d), we have
$$S^{k+1}_\z=S^1_{\z^0+p\r}(S^1_{\z^1+p\r})^{(1)}\do(S^1_{\z^{k-1}+p\r})^{(k-1)}(S^1_{\z^k})^{(k)}.\tag a$$

\subhead 6 \endsubhead
According to \cite{\ANI, 5.2(a)}, for $\nu\in X^+_{red}$ we have
$$S_{\nu+p\r}=S^1_{\nu+p\r}\text{ provided that }p\gg0.\tag a$$
It is likely that for $\nu\in X^+_{red}\cap X^{++}$ we have
$$S_\nu=S^1_\nu\text{ provided that }p\gg0.\tag b$$
Note that (b) is a very special case of Conjecture 5.1 in \cite{\ANI}.

\proclaim{Proposition 7}Assume that $k\ge0$ and $\z\in X^+_k$. Assume that $p\gg0$ and that 6(b) holds. Then
$S_\z=S^{k+1}_\z.$
\endproclaim
If $k=0$ this is just the assumption 6(b). Assume now that $k\ge1$.
Using 3(c) and 5(a) we see that it is enough to show
$$\align&S_{\z^0+p\r}S_{\z^1+p\r}^{(1)}\do S_{\z^{k-1}+p\r}^{(k-1)}S_{\z^k}^{(k)}\\&
=S^1_{\z^0+p\r}(S^1_{\z^1+p\r})^{(1)}\do(S^1_{\z^{k-1}+p\r})^{(k-1)}(S^1_{\z^k})^{(k)}.\tag a\endalign$$
This follows from the equalities $S_{\z^j+p\r}=S^1_{\z^j+p\r}$ for $j=0,1,\do,k-1$ (see 6(a)) and
$S_{\z^k}=S^1_{\z^k}$ which holds by the assumption 6(b).

\proclaim{Corollary 8}In the setup of Proposition 7 we have
$$S_\z=\sum_{\mu_0,\mu_1,\do,\mu_k\in X^{++}}
x_{\mu_k,\z^k}x_{\mu_{k-1},\z^{k-1}+p\mu_k}\do x_{\mu_1,\z^1+p\mu_2}x_{\mu_0,\z^0+p\mu_1}S^0_{\mu_0}.$$
\endproclaim
Using 4(a) repeatedly, we have
$$S^1_\z=\sum_{\mu_0\in X^{++}}x_{\mu_0,\z^0+p\z^1+p^2\z^2+\do}S^0_{\mu_0},$$
$$S^2_\z=\sum_{\mu_1\in X^{++}}x_{\mu_1,\z^1+p\z^2+\do}S^1_{\z^0+p\mu_1}
=\sum_{\mu_0,\mu_1\in X^{++}}x_{\mu_1,\z^1+p\z^2+\do}x_{\mu_0,\z^0+p\mu_1}S^0_{\mu_0},$$
$$\align&S^3_\z=\sum_{\mu_2\in X^{++}}x_{\mu_2,\z^2+p\z^3+\do}S^2_{\z^0+p\z^1+p^2\mu_2}\\&
=\sum_{\mu_0,\mu_1,\mu_2\in X^{++}}x_{\mu_2,\z^2+p\z^3+\do}x_{\mu_1,\z^1+p\mu_2}x_{\mu_0,\z^0+p\mu_1}
S^0_{\mu_0}.\endalign$$
Continuing we get
$$S^{k+1}_\z=\sum_{\mu_0,\mu_1,\do,\mu_k\in X^{++}}
x_{\mu_k,\z^k}x_{\mu_{k-1},\z^{k-1}+p\mu_k}\do x_{\mu_1,\z^1+p\mu_2}x_{\mu_0,\z^0+p\mu_1}S^0_{\mu_0}.$$
It remains to use Proposition 7.

\subhead 9 \endsubhead
For any $\z\in X^{++}$ we can write $S_\z=\sum_{\mu\in X^{++}}r_{\mu,\z}S^1_\mu$ where $r_{\mu,\z}\in\NN$. 

\proclaim{Corollary 10}In the setup of Proposition 7 assume that $\mu\in X^{++}$ satisfies $r_{\mu,\z}\ne0$.
Then $\mu=\z\mod pX$.
\endproclaim
Using Corollary 8 and 4(a) we see that
$$r_{\mu,\z}=\sum_{\mu_1,\do,\mu_k\in X^{++};\mu=\z^0+p\mu_1}
x_{\mu_k,\z^k}x_{\mu_{k-1},\z^{k-1}+p\mu_k}\do x_{\mu_1,\z^1+p\mu_2}.$$
It follows that $\mu=\z^0+p\mu_1$ for some $\mu_1\in X^{++}$. Since $\z=\z^0\mod pX$, it follows that 
$\mu=\z\mod pX$. The corollary is proved.

\subhead 11 \endsubhead
In view of Corollary 10, one could hope that for any $\z,\mu\in X^{++}$ such that $r_{\mu,\z}\ne0$ we have
$\mu=\z\mod pX$. Unfortunately, this is contradicted by example (i) in \cite{\JE}.

\enddocument